% Logic Eprints
%Submitted 0738 Fri Aug 06, 1993 by: goldstrn@math.fu-berlin.de (martin goldstern)
%logic/goldstern/shoenfield.tex
%

\magnification\magstep1
\font\bigfont=cmr17
\font\small=cmsl10 at 10 true pt	

\advance\topskip \baselineskip

\parindent0cm
\def\begindent{\par\begingroup\par\advance\parindent1cm\icount0 }
\def\endent{\par\endgroup\par}
\def\ite #1 {\item{(#1)}}
\def\itm{\advance\icount1 \ite{\number\icount} }
\newcount\icount

\newcount\thmcount
\def\prop#1:{\medskip\advance\thmcount1 
\noindent{\bf \number\thmcount. #1:} \ignorespaces} 

\def\lprop#1 #2:{\prop#2:\label{#1}\ignorespaces}

\def\label#1{\expandafter\edef\csname LABEL#1\endcsname{\number\thmcount}}
\def\ref#1{\csname LABEL#1\endcsname}

\def\<#1>{\hbox{$\langle #1\rangle$}}

\def\sq{{}^\omega \omega }
\def\two#1{{}^{#1}2}
\def\C{{\bf C}}
\let\cut\cap
\let\bigcut\bigcap

\def\div{\omega {\nearrow} \omega}
\def\B{\hbox{\rm I}\!\hbox{\rm B}}
\def\0{{\bf 0}}
\def\1{{\bf 1}}

\abovedisplayskip0.5\abovedisplayskip
\belowdisplayskip0.5\belowdisplayskip

\newcount\refcount

\def\key#1{\advance\refcount1 
\edef#1{\number\refcount}}

\key\DrWi
\key\fkt
\key\grabner
\key\grill
\key\multiple
\key\lusin
\key\mosch
\key\neumann
\key\shoen
\key\randomsurvey
\key\winkler

^^L

\medskip
{{\bigfont
\obeylines\everypar{\hskip\parfillskip}
\parskip20pt

\vglue 3cm

An application of  Shoenfield's absoluteness theorem
to the theory of uniform distribution
\bigskip\bigskip

\bigfont
 Martin Goldstern
\bigskip
{\rm October 1992/January 1993}

}

\vfill

{\leftskip0.15\hsize\rightskip0.15\hsize\noindent ABSTRACT.  Theorem:
If $\C$ is a Polish probability space, $W \subseteq \omega ^\omega
\times \C$  a Borel set  whose sections $W_x$ ($x\in
\omega^\omega$) have measure one and  are decreasing ($x\le x' \to
W_x \supseteq W_{x'}$), then the set $\bigcut_x W_x$ has measure one. 

We give two proofs of this theorem: 
 one  in the language of set theory,
the other in the language of probability theory,   and we apply the
theorem to a question on completely uniformly distributed sequences.

}

\pageno=0\nopagenumbers
\bigskip\bigskip
Supported by DFG grant Ko 490/7-1
\eject

}

\headline{\small\hfil
Shoenfield's absoluteness theorem and uniform distribution\hfil}

\beginsection{Introduction.}

The law of large numbers implies that almost every infinite sequence
in $\{0,1\} $ is uniformly distributed, i.e., for almost all
 $t\in \two \omega$ we have $\lim_{n\to \infty} D(t,n)=0$, where
$2=\{0,1\}$, $\omega = \{0, 1, 2, \ldots\}$, and the ``discrepancy''
$D$ is defined by
 $D(t,n) = 2\cdot\left| { |\{i<n:t(i)=0\}|\over n}-{1\over 2}\right|$,
and ``almost all'' refers to the product measure on $\two\omega$
induced by the equidistributed measure  on $\{0,1\}$, i.e.,
$\mu(\{0\})=\mu(\{1\})={1\over 2}$. 

A similar argument shows that for all $k$, almost all sequences are
$k$-uniformly distributed, i.e., for all $k$, for almost all $t\in
\two\omega$ we have $\lim_{n\to \infty} D_k(t,n)=0$, where 
$$D_k(t,n):= 
2 ^k   \cdot \max_{w\in \two k}	\left| 
     { |\{i< n: (t(i), \ldots, t(i+k-1)) = w \} | 
      \over 
     n}
- {1 \over 2^k} 
\right| $$
Since the measure is countably additive, we can switch the quantifiers
``for all $k$'' and ``for almost all $t$'' and get: 
$$\hbox{ for almost all $t\in \two \omega $, for all $k$, $t$ is
$k$-uniformly distributed}$$

If we replace the constant $k$ by a function $s\in \sq$ we are led to
the following concept:  A sequence $t\in \two\omega$ is called
$s$-uniformly distributed, if $\lim_{n\to \infty}D_{s(n)}(t,n)=0$.

Let $R_s$ be the set of $s$-uniformly distributed sequences. It is
easy to see that $s$ must not grow faster than $\log n$ (in this
paper, $\log$ always means the logarithm to the base 2)
if the set $R_s$ should be nonempty.  Flajolet-Kirschenhofer-Tichy and
Grill determined
 the exact  bound for $s$  below which still almost all sequences are
$s$-uniformly distributed: 

\prop Theorem ([\fkt], [\grill]): $\mu(R_s) = 1 $ 
iff $\varphi_s\to \infty$,
 where $\varphi_s(n) = \lfloor \log n - \log\log n - s(n)\rfloor$. 

See also [\grabner] for a generalization of this result. 
\bigskip

Now let $$ R = \bigcut _{s} R_s$$
where $s$ ranges over all functions with $\varphi_s\to \infty$. Since
$R$ is an uncountable intersection of measure one sets, it is a priori
not clear if $R$ has measure one, or indeed if $R$ can be nonempty.

We will show below that $R$ must have measure $1$.  That is, as in the
remark about $k$-uniformly distributed, we can again switch
quantifiers and get

\prop Theorem: For almost all $t\in \two \omega $, for all $s$ for
which $\varphi_s\to \infty$ we have: $t$ is $s$-uniformly distributed.

%\smallskip
%The same is of course true for  all spaces ${}^\omega m$, $m\ge 2$,
%if we read $\log$ as $\log_m$. 

\medskip
We will give two proofs of this theorem:  One for logicians, using
the technique of forcing and Shoenfield's absoluteness theorem for
$\bf\Sigma^1_2 $-sets, and
the second for those mathematicians who are more familiar with the
language of probability theory.  Here we use von Neumann's theorem
that Borel  sets can be uniformized by measurable functions.

\lprop notation Notation: 
 $\omega$ is the set of natural numbers.  $2 = \{0,1 \}$. 

$\sq $ is the Baire space ( = all functions from $\omega$
to $\omega$) with the usual metric and topology.  The variable $x$
always ranges over elements of $\sq$.  This spaces is naturally
ordered by $x_1 \le x_2 \ \Leftrightarrow\  \forall n \, x_1(n) \le x_2(n)$. 

$\div $ is the set of  sequences which diverge to
infinity, a subspace of $\sq$. The variable $y$ always ranges over
elements of $\div$. 

%$\C$ (where $m = \{0,\ldots, m-1\}$)  is the $m$-Cantor space.  With
%the obvious metric it is a compact metric space, and it carries the
%usual probability measure (the product measure induced by the
%equidistributed measure on $m$). 

Let $\C$ be any perfect Polish space (i.e., complete metric
separable) with a complete probability measure $\mu$ (i.e., $\mu$ is
${\sigma}$-additive, all Borel sets are measurable, and all subsets of
measure zero sets are measurable).   The reader may think of
$\C=\two\omega$, the set of all functions from $\omega$ to $2$.  The
variable $t$ always ranges over $\C$.   

Analytic sets (also called  $\bf {\Sigma}^1_1$ sets) are those subsets
of $\C$ which can be obtained as projections of Borel
sets in $\C \times  \sq$ (or equivalently, are a continuous image of
$\omega^ \omega$). A subset of $\C$ is coanalytic or
$\bf\Pi^1_1$  iff its complement is analytic.  (All propositions below
will remain true if ``$\bf {\Sigma}^1_1$'' and ``$\bf\Pi^1_1$'' are
replaced by ``Borel''.) 
It is well known that all analytic sets are measurable. ([\lusin], see
also [\mosch, 2H.8])

All facts from descriptive set theory that we use here can be found in
Moschovakis book [\mosch].  For facts and references about forcing
consult [\multiple].

\prop Acknowledgement:  Thanks to Reinhard Winkler for interesting
dicussions about random numbers.

\bigskip

\beginsection{The main lemma and its application}

\lprop mainlemma Main Lemma:  Assume that for 
each $y\in\div$ we have a  set $U_y
\subseteq \C$ such that  
\begindent
\ite A The set $\{(y, t): y\in\div, t\in U_y\}$ is a   $\bf\Pi^1_1$
  subset of $\subseteq  \div \times \C$. 
\ite B For all $y\in \div$, $\mu(U_y) = 1$
\ite C If $y \le y'$ (that means $\forall  n\, y(n) \le
       y'(n)$), then $U_y \subseteq U_{y'}$
\endent
Then $\mu(\bigcut_{y\in \div } U_y) =1$. 

%\prop Note: There are
%       obvious generalizations to ${\sigma}$-finite measures on Polish
%       spaces: We can replace ``$=1$'' by ``$>0$'', or also by
%       $=\infty$''. 
%

\medskip

We will actually not prove the lemma itself, but an equivalent
version: 

\lprop equiv Lemma:  Assume that for
each $x\in \sq$ we have a  set $W_x \subseteq \C$ such
that  \begindent
\ite A' The set $\{(x, t): x\in\sq, t\in W_x\}$ is a $\bf \Pi^1_1$ set 
  $\subseteq  \sq  \times \C$. 
\ite B' For all $x\in \sq$, $\mu(W_x) = 1$
\ite C' If $x \le x'$, then $W_x \supseteq W_{x'}$ 
\endent
Then $\mu(\bigcut_{x\in \sq} W_x) =1$. 

\medskip
Proof that lemma \ref{equiv} implies lemma \ref{mainlemma}: 
 For every unbounded  $x\in \sq$ we define its ``inverse'' function
$x^*\in \div$ by 
$$x^*(n)= \min\{ m: x(m) \ge n\}$$
Now given a family $(U_y: y\in \div)$ let $W_x := U_{x^*}$ for all
unbounded $x\in \sq$, and $W_x=\C$ otherwise.  

We claim $\bigcut_y U_y = \bigcut_x W_x$.  The inclusion
``$\subseteq$'' is clear.  For the converse inclusion, it is enough to
see 
$$ \forall y \in \div \, \exists x \in \sq: x ^* \le y$$
Given $y\in \div$, we can find a sequence $x$ which increases so fast
that for all $n$, $x(y(n)) \ge n$, e.g., $x(m) := \max \{ k: y(k) \le
m\}$.  So $x^*(n)  = \min\{ i: x(i) \ge n \} \le y(n)$. 

Hence  $\bigcut _x W_x = \bigcut _y U_y$ has measure 1.

\lprop app Application:  let $\C = \two \omega$ with the usual product
measure.  
Let $s$ range over functions in $\sq$
satisfying $ \varphi_s \ge0 $ and  $\varphi_s\to \infty$, where
 $\varphi_s=\lfloor \log n - \log\log n - s(n) \rfloor$.  Let $R_s$ be
the set of $s$-uniformly distributed sequences in $2 ^\omega$, and let
$R:=\bigcut_s R_s$.   It is easy to check that $s \le s'$ implies
$R_{s'} \subseteq R_s$, and $\varphi_{\varphi_s} = s$. 

For each function $y\in \div$ let 
$$U_y = \cases{ R_{\varphi_y}& if $y \le \log -\log\log $\cr
               2 ^\omega         & otherwise\cr}$$
Now check that all the assumptions of the main lemma are satisfied.
If $\varphi_s\to \infty$, then letting $y:= \varphi_s  $ we have
$s= \varphi_y$, so $R_s\supseteq R_{\varphi_y} = U_y$.  Hence $R =
\bigcut_s R_s$ has measure 1.

\prop Remark:  Note that the inclusion in lemma \ref{mainlemma}(C) or
lemma \ref{equiv}(C') 
cannot be replaced by an ``almost'' inclusion (modulo measure zero
sets or even modulo finite sets): Let $F$ be a Borel isomorphism
between $\sq$ and $\C$ (see [\mosch, 1G4]), and let 
$W_x = \C - \{F(x)\}$, then $\bigcut_x W_x =\emptyset$.

\prop Other applications:  Almost every metric result in which a 
$o(\,\, )$ appears can be strenghtened by using lemma~\ref{mainlemma}. 
   For example,  Drmota-Winkler
showed that \def\oo{\sqrt{ N / \log N}}
$$\hbox{ If $s(N) = o(\oo)$, then $\mu(T_s)=1$,} $$
where $T_s$  is the set of $s$-uniformly distributed sequences  
$x\in [0,1]^\omega$ (see [\DrWi] for definitions and references), and
$s(N) = o(\oo)$ means that $\lim_{N} s(N) \left/\oo\right. = 0$.  

Now for $y\in \div$ let $s_y(N)= \left\lfloor \left. \oo \right/
y(N)\right\rfloor$, and let 
$U_y = T_{s_y}$. Applying the lemma we   immediately get
that 
$$\mu\bigl(\bigcut_s T_s \bigr) = 1 $$
where the intersection is taken over all $s$ satisfying $s(N) =
o(\oo)$.

\beginsection{Proof of lemma \ref{equiv} for logicians}

The set $ W := \{ t: \forall x \, t\in
W_x\}$ is a $\bf\Pi^1_1$ 
set, hence measurable.  Assume that $\mu(W)< 1$.   So there exists a
Borel ($G_{\delta}$) set $A$ with $\mu(A)<1 $ such that $W \subseteq
A$.  Note that the statement
$$ (*) \qquad\qquad 
 \forall t:  \left[ \forall x \, t\in W_x \ \Rightarrow \ t\in A\right]$$
is a $\bf\Pi^1_2$-statement, hence absolute between any two transitive
universe with the same ordinals ([\shoen], see also [\mosch, 8F10]). 

Now let $r$ be a random real over the universe $V$, $r\in \C - A$, 
 and work in $V[r]$.  Random forcing is $\sq$-bounding, i.e., we have 
$$ \forall x\in \sq\, \exists x'\in \sq\cut V : x \le x'$$
Since the condition (C') is also absolute, this implies $ \forall x\in
\sq\, \exists x'\in \sq\cut V: W_x \supseteq W_{ x'}$, and so 
$$  \bigcut _{x\in \,\sq} W_x = 
\bigcut_{x'\in \,\sq\cut V} W_{x'}$$
Moreover, 
$$ \forall x'\in \sq\cut V : r \in W_{x'}$$
(since $r$ is a random real and each $W_{x'}$ contains a Borel set of
  full measure), so   together we get  
$$ \forall x\in \sq :  r \in W_x. $$
Since $(*)$ holds also in $V[r]$, we get $V[r] \models r\in A$, a
contradiction.

\beginsection{Proof of lemma \ref{equiv} for  probabilists}

Now we repeat our argument in the language of probability theory.  
Our main tool is  von Neumann's selection theorem
([\neumann], see also [\mosch, 4E9]).

 We let our probability space be $\C$. 

We will consider random variables (= measurable functions)  $T$ from
$\C$ into $\C$, and random variables $X$ from $\C$ into $\sq$.  

\prop Definition:  Let $B \subseteq \C\times \sq$.  We say that $f$
``uniformizes'' $B$ iff
\begindent
\ite 1 $f \subseteq B$
\ite 2 $f$ is a partial function. 
\ite 3 $dom(f) = \{t: \exists x\, (t,x)\in B\}$. 
\endent

The von Neumann selection theorem says  that every $\bf {\Sigma}^1_1$
set can  be uniformized by a measurable 
function.    (Note
that the possibly better known $\bf\Pi^1_1$ uniformization theorem will
in general not yield measurable functions [\mosch, 5A7].) 

 Hence we get:

\lprop uniform Lemma:  
 Let $B \subseteq \C \times \sq$ be a $\bf \Sigma^1_1 $ set.
Assume that $\forall t\in \C\, \exists x\in \sq \, (t,x)\in B$. 
Then 
\begindent
\ite 1 $\exists X \mu [(I,X)\in B] =1 $, where $I$ is the identity
      function from $\C$ to $\C$, and $X$ ranges over all random
      variables in $\sq$. 
\ite 2 Moreover, $\forall T \,\exists X  \mu[(T,X)\in B] = 1 $.   Here,  $T$
     and $X$ range over random variables in $\C$ and $\sq$, respectively. 
\endent

Proof: (1) ---  There is a measurable function $f \subseteq B$ with domain
      $\C$.  So for all $ t\in \C$, $(t,f(t))\in B$.  Let $X:=f$. 

(2) --- Let $f$ be the function from (1), and let $X = f \circ T$.

\prop Remark: Those logicians who have not skipped this section
will notice that (2) is again Shoenfield's theorem, since random
variables naturally correspond to $\B$-names of reals (where $\B$ is
the random algebra).

		\bigskip
We will also need the following easy fact about random variables
$X:\C\to \sq$ (which corresponds to the fact that random forcing is
$\sq$-bounding):

\prop Fact: Let $X:\C\to \sq$ be a random variable.  Then there is  a
  family $(x_n:n\in \omega)$ of functions in $\sq$ such that 

$$ \mu  [ \exists n \,  X \le x_n ]  = 1 $$

Proof: Call a Borel set $A \subseteq \C$ of positive measure ``good'' if
there is a function $x_A\in \sq$ such that $A \subseteq [X\le x_A]$.
If there are finitely many good sets covering $\C$ (up to a measure
zero set) then we are done.
Otherwise, let $(A_n:n\in \omega)$ be a maximal antichain of good
sets, i.e., a maximal family of sets which are 
\begindent
\ite 1 good
\ite 2 pairwise disjoint
\endent
Such a sequence can be found using Zorn's lemma. 
Again, if $\bigcup_n  A_n$ covers $\C$ up to a measure zero set
then we are done.  So assume that $A:=\C - \bigcup _n A_n$ has
positive measure $\varepsilon$.  Define a sequence $(k_n:n\in \omega)$
of natural numbers such that for each $n$ the set $[X(n)> k_n]$ has
measure $< \varepsilon/2^{n+3}$. Thus, the set
 $[\forall n\, X(n)\le k_n] \cap A$ has measure $\ge \varepsilon -
\varepsilon/4 >0$, contradicting the maximality of the family
$(A_n:n\in \omega)$.

\prop Second proof of lemma~\ref{equiv}:
Let $I:\C \to \C$ be the identity function.  Clearly, for any $x\in
\sq$ we have $\mu [I \in W_x] =\mu (W_x) =  1$.

 We now claim that the same is still true if we replace $x$ by a
random variable $X$.   Indeed, let $X$ be a random variable, $X:\C \to
\sq$.  Then we can find a family of  functions $(x_n:n\in \omega)$ such
that
$$\mu [ \exists n \,  X \le x_n]=1  $$
Using the anti-monotonicity of the family $(W_x:x\in \sq)$ we get 
 $$\mu [\exists n \, W_X \supseteq W_{x_n} ]=1 $$
and so: $\mu [W_X \supseteq \bigcut_n W_{x_n}] = 1$.  Hence 
$$ (*) \qquad \qquad 
\mu [I \in W_X ] \ge \mu  [ I\in   \bigcut_n W_{x_n}]=
\mu\bigl( \bigcut_n W_{x_n}\bigr) = 1$$ 
 for all random variables $X$.  So for all $X$ $\mu[I\notin W_X]=0$,
      hence by the contrapositive of lemma \ref{uniform} it is
      impossible that $\forall t \,\exists x\,\, t\notin W_x$. Therefore 
$$ \exists t\, \forall x \,\, t \in W_x$$
So the set $\bigcut_x W_x$ is nonempty.  Relativizing the above
argument to any positive Borel set $A$ we can show that the set
$A\cut \bigcut_x W_x$ is nonempty.  
Hence $W$ has outer measure 1.  As remarked above, $W$ must be
measurable, so indeed $\mu(W)=1$.

\bigskip\bigskip\bigskip\bigskip\bigskip

\beginsection References.

\baselineskip\normalbaselineskip
\def\key#1{\smallskip
\hangindent0.5cm\hangafter1 [#1] }

\key\DrWi M.~Drmota and R.~Winkler, $s(N)$-uniform distribution modulo 1.

\key\fkt P.~Flajolet, P.~Kirschenhofer and R.F.~Tichy, Deviations from
	uniformity in 
	random strings, Probab.\ Th.\ Rel.\ Fields 80 (1988), 139--150. 

\key\grabner P.~Grabner, Block Distribution in Random Strings, Annales
	de l'Institut Fourier. 

\key\grill K.~Grill, A note on Randomness, Stat.\ and Probab.\ Letters. 

\key\multiple T.~Jech, Multiple forcing,  Cambridge Tracts in Mathematics
88, Cambridge 1986.

\key\lusin N.~Lusin, Sur la classification de M.~Baire,
       Comptes Rendus Acad.~Sci.~Paris 164 (1917) 91--94.

\key\mosch Y.N.~Moschovakis, Descriptive Set Theory, Studies in Logic and
       the foundation of Mathematics, vol 100, North-Holland 1980. 

\key\neumann J.~von Neumann, On rings of operators. Reduction theory.
Ann.\ Math.\ 50 (1949) 448--451.

\key\shoen J.R.~Shoenfield, The problem of predicativity, Essay on the
       foundation of Mathematics, Magnes Press, Hebrew Univ. Jerusalem
       (1961) 132--139.  

\key\randomsurvey R.F.~Tichy, 
Zur Analyse und Anwendung von Zufallszahlen, Jahrbuch d.\ Kurt Goedel
	Ges.\ 109--116 (1990

\key\winkler R.~Winkler, Some remarks on Pseudorandom Sequences. 

\bigskip\vfill

\obeylines\frenchspacing\baselineskip\normalbaselineskip\parskip0pt
\vbox{
Martin Goldstern
2. Mathematisches Institut
FU Berlin
Arnimallee 3
1000 Berlin 33
}
\smallskip
\tt goldstrn@math.fu-berlin.de
\eject

\bye